\documentclass[11pt]{amsart}

\usepackage[T1]{fontenc}
\usepackage{amsmath,amssymb,amsthm}
\usepackage{enumitem}
\usepackage[colorlinks=true,
            linkcolor=blue,
            citecolor=blue,
            urlcolor=blue]{hyperref}

\newtheorem{theorem}{Theorem}[section]
\newtheorem{lemma}[theorem]{Lemma}
\newtheorem{proposition}[theorem]{Proposition}
\newtheorem{corollary}[theorem]{Corollary}

\theoremstyle{definition}
\newtheorem{definition}[theorem]{Definition}
\newtheorem{problem}[theorem]{Problem}

\theoremstyle{remark}
\newtheorem{remark}[theorem]{Remark}
\newtheorem{example}[theorem]{Example}

\DeclareMathOperator{\coz}{coz}
\DeclareMathOperator{\Min}{Min}
\DeclareMathOperator{\Z}{Z}
\DeclareMathOperator{\K}{K}

\newcommand{\CZ}{\mathcal Z}

\newcommand{\F}{\mathcal F}

\title[$n$-absorbing ideals of $C(X)$: minimal primes and common zero sets]
{$n$-absorbing ideals of rings of continuous functions I:
minimal primes, common zero sets, and non-$z$ ideals}

\author{Ali Akbar Estaji}
\address{Faculty of Mathematics and Computer Science,
Hakim Sabzevari University, Sabzevar, Razavi Khorasan Province,
Postal Code 9617976487, Iran}


\date{}

\begin{document}

\begin{abstract}
Let $C(X)$ be the ring of real-valued continuous functions on a
Tychonoff space $X$, and for an ideal $I$ let $\K(I)$ be the set of
common zeros of its members.
We show that a proper $z$-ideal of $C(X)$ is $n$-absorbing if and only
if it has at most $n$ minimal primes, if and only if it is an
intersection of at most $n$ pairwise incomparable primes, if and only if
its filter of zero sets is $n$-absorbing; the algebraic part of this
equivalence is valid for radical ideals of any commutative ring.
We prove that $|\K(I)|\le n$ for every $n$-absorbing ideal $I$ and that
this bound is sharp.
We further prove that the $z$-closure of every $n$-absorbing ideal
is an $n$-absorbing $z$-ideal and has the same common zero set.
Finally, every $n$-absorbing ideal of $C(X)$ is a $z$-ideal if and
only if $X$ is a $P$-space. Moreover, for every non-unit
$f\in C(X)$, the principal ideal $(f)$ is $n$-absorbing if and only
if $\Z(f)$ is clopen and $(f)$ has at most $n$ minimal prime ideals.
\end{abstract}

\maketitle

\section{Introduction}

The theory of $n$-absorbing ideals provides a natural extension of prime
ideal theory. A proper ideal $I$ of a commutative ring $R$ is called
\emph{$n$-absorbing} if, whenever
$a_1\cdots a_{n+1}\in I$, the product of some $n$ of the elements
$a_1,\ldots,a_{n+1}$ belongs to $I$. The case $n=2$ was introduced by
Badawi~\cite{Badawi}, and the general notion was developed by Anderson
and Badawi~\cite{AndersonBadawi}. Among the basic structural results,
the intersection of $n$ prime ideals is $n$-absorbing, and every
$n$-absorbing ideal has at most $n$ minimal prime ideals
\cite[Theorems~2.1(c) and~2.5]{AndersonBadawi}. For $n=2$, Badawi also
showed that the radical of a $2$-absorbing ideal is either prime or the
intersection of two minimal prime ideals~\cite{Badawi}.

For a Tychonoff space $X$, the ring $C(X)$ of real-valued continuous
functions carries additional structure that is invisible in the
abstract ring-theoretic definition. In particular, $z$-ideals are
controlled by zero sets and correspond to filters in the lattice of
zero sets. Bilgin, Demir, and Oral established the correspondence
between $2$-absorbing $z$-filters and $2$-absorbing $z$-ideals of
$C(X)$~\cite{BilginDemirOral}. Mason's work supplies another essential
ingredient: prime ideals minimal over a $z$-ideal are themselves
$z$-ideals~\cite[Theorem~1.1]{Mason}. The topology of minimal-prime
spectra has also been studied in the general ring-theoretic setting
by Henriksen and Jerison~\cite{HenriksenJerison}. Related absorbing
conditions for lattice ideals and filters have been investigated by
Ebrahimi Atani and Sedghi Shanbeh Bazari~\cite{EbrahimiAtaniSedghi}
and by Estaji and Haghdadi~\cite{EstajiHaghdadi}. These works provide
a broader context for the present study, which focuses on the
interaction between the $n$-absorbing property and the topology
encoded by $C(X)$.

The central object in this paper is the common zero set
$
  \K(I)=\bigcap_{f\in I}\Z(f),
$
which records the points at which every member of an ideal $I$ vanishes.
For a prime ideal $P$ of $C(X)$, this set is either empty or a singleton.
For an $n$-absorbing ideal, the number of minimal primes is bounded by
$n$; however, several distinct minimal prime $z$-ideals may lie over
the same point of $X$, so the common zero set does not in general
record the full minimal-prime structure. This finer structure is the
subject of the companion paper~\cite{EstajiB}, which builds on the
results proved here.

Our first results establish the precise algebraic and zero-set
framework. We extend the correspondence of Bilgin, Demir, and Oral from
$2$-absorbing to arbitrary $n$-absorbing $z$-ideals, and record the
radical-ideal characterization in terms of minimal primes. The bound
on the number of minimal primes is already known for arbitrary
$n$-absorbing ideals by Anderson and Badawi; the localization argument
included here gives a self-contained proof in the radical setting.
We then prove a topological bound, namely
$|\K(I)|\leq n$ for every $n$-absorbing ideal $I$ of $C(X)$, and show
that the bound can be attained whenever $X$ has at least $n$ points.
We also show that the $z$-closure of an arbitrary $n$-absorbing ideal
is again $n$-absorbing, using its minimal-prime decomposition.
Finally, we consider ideals that are not necessarily $z$-ideals,
recall the bound on the number of minimal primes of their radicals,
characterize the principal $n$-absorbing ideals, and record the
resulting characterization of $P$-spaces.

The paper is organized as follows. Section~2 collects the required
facts about zero sets, $z$-ideals, and minimal primes. Section~3 relates
$n$-absorption to zero-set filters and minimal-prime decompositions.
Section~4 studies the common zero set. Section~5 treats arbitrary
$n$-absorbing ideals, and Section~6 lists further questions.

\section{Preliminaries}

All rings are commutative with identity, and all ideals are proper
unless stated otherwise. We collect standard facts about $z$-ideals and
refer to the literature for proofs. For $f\in C(X)$ put
$\Z(f)=\{x\in X: f(x)=0\}$ and $\coz(f)=X\setminus\Z(f)$.
Since $C(X)\cong C(Y)$ for a Tychonoff space $Y$ depending on $X$
(see~\cite[Chapter~3]{GillmanJerison}), we assume throughout that
\emph{$X$ is a Tychonoff space} (completely regular and Hausdorff).
The family of zero sets $\CZ(X)$ is a bounded distributive lattice
under $\cap$ and $\cup$, because $\Z(f)\cap\Z(g)=\Z(f^2+g^2)$ and
$\Z(f)\cup\Z(g)=\Z(fg)$.

For $x\in X$ let $M_x=\{f\in C(X):f(x)=0\}$; it is a maximal ideal and
a $z$-ideal. For an ideal $I$ put
$\K(I)=\bigcap_{f\in I}\Z(f)$ and $\CZ[I]=\{\Z(f):f\in I\}$.
An ideal $I$ is a \emph{$z$-ideal} if $f\in I$ and
$\Z(f)\subseteq\Z(g)$ imply $g\in I$.
A \emph{filter} of $\CZ(X)$ is a nonempty $\F\subseteq\CZ(X)$ with
$\varnothing\notin\F$, closed under finite intersections and upward
closed in $\CZ(X)$; it is \emph{prime} if $A\cup B\in\F$ implies
$A\in\F$ or $B\in\F$.
If $X$ is compact, every proper ideal of $C(X)$ is fixed, i.e.\
$\K(I)\neq\varnothing$ \cite[Chapter~4]{GillmanJerison}.
For an ideal $I$, $\Min(I)$ denotes the set of prime ideals minimal over
$I$.

The \emph{$z$-closure} of an ideal $I$ is the smallest $z$-ideal containing $I$:
\[
I_z=\bigl\{f\in C(X):\ \exists\,g\in I\text{ with }\Z(g)\subseteq\Z(f)\bigr\}.
\]
Equivalently, $I_z=\mathrm{Z}^{-1}[\mathrm{Z}[I]]$. Thus $I$ is a $z$-ideal
if and only if $I=I_z$.

\begin{lemma}
\label{lem:basic}
\begin{enumerate}[label=\textup{(\alph*)}]
\item Every $z$-ideal of $C(X)$ is radical, and an intersection of
$z$-ideals is a $z$-ideal.
\item If $I$ is a $z$-ideal, then $f\in I\iff\Z(f)\in\CZ[I]$, and
$\CZ[I]$ is a filter of $\CZ(X)$.
\item A $z$-ideal $I$ is prime iff $\CZ[I]$ is a prime filter.
\end{enumerate}
\end{lemma}

\begin{proof}
Statement \textup{(a)} is immediate from the definition. For
\textup{(b)} and \textup{(c)} see \cite[Chapters~2 and~14]{GillmanJerison}.
\end{proof}

\begin{lemma}
\label{lem:elementary}
Let $R$ be a commutative ring.
\begin{enumerate}[label=\textup{(\alph*)}]
\item The intersection of $n$ prime ideals, not necessarily distinct, is $n$-absorbing.
\item If $P_1,\dots,P_m$ are pairwise incomparable prime ideals, then
$\Min(P_1\cap\dots\cap P_m)=\{P_1,\dots,P_m\}$; in particular, for
$m\ge2$ the intersection is not prime.
\end{enumerate}
\end{lemma}

\begin{proof}
For \textup{(a)} see \cite{AndersonBadawi}; the case $n=2$ is
\cite[Theorem~2.3]{Badawi}. For \textup{(b)} see
\cite[Chapter~1]{AtiyahMacdonald}.
\end{proof}

\begin{lemma}[Mason]
\label{lem:min-z}
Let $I$ be a radical $z$-ideal of $C(X)$. Every prime ideal minimal over
$I$ is a $z$-ideal. Moreover, if $I$ is proper then $\Min(I)\neq\varnothing$
and $I=\bigcap_{P\in\Min(I)}P$.
\end{lemma}

\begin{proof}
The first statement is \cite[Theorem~1.1]{Mason}. The second is the
standard fact that a radical ideal equals the intersection of the
primes containing it, and every prime containing $I$ contains a minimal
one; see \cite[Chapter~1]{AtiyahMacdonald}.
\end{proof}

\begin{definition}
\label{def:P-point}
A point $x\in X$ is a \emph{$P$-point} if every zero set containing $x$
is a neighbourhood of $x$; $X$ is a \emph{$P$-space} if all its points
are $P$-points.
\end{definition}

\section{\texorpdfstring{$n$-absorbing $z$-ideals}{n-absorbing z-ideals}: zero sets and minimal primes}

\begin{definition}
Let $n\ge1$. A filter $\F$ of $\CZ(X)$ is \emph{$n$-absorbing} if, for
all $A_1,\dots,A_{n+1}\in\CZ(X)$, $A_1\cup\dots\cup A_{n+1}\in\F$
implies that the union of some $n$ of the sets $A_i$ lies in $\F$.
For $n=1$ these are the prime filters, and for $n=2$ the
$2$-absorbing $z$-filters of~\cite{BilginDemirOral}.
\end{definition}

\begin{proposition}[{compare~\cite{BilginDemirOral} for $n=2$}]
\label{prop:zero-set}
A proper $z$-ideal $I$ of $C(X)$ is $n$-absorbing if and only if
$\CZ[I]$ is an $n$-absorbing filter.
\end{proposition}

\begin{proof}
Assume $I$ is $n$-absorbing and
$A_i=\Z(f_i)$ with $A_1\cup\dots\cup A_{n+1}=\Z(f_1\cdots f_{n+1})\in\CZ[I]$.
By Lemma~\ref{lem:basic}(b), $f_1\cdots f_{n+1}\in I$, so for some $k$
the product $\prod_{l\ne k}f_l$ lies in $I$, and its zero set
$\bigcup_{l\ne k}A_l$ lies in $\CZ[I]$.
Conversely, assume $\CZ[I]$ is $n$-absorbing and
$f_1\cdots f_{n+1}\in I$. Then $\bigcup_i\Z(f_i)\in\CZ[I]$, so for some
$k$ we have $\Z\bigl(\prod_{l\ne k}f_l\bigr)=\bigcup_{l\ne k}\Z(f_l)\in\CZ[I]$,
and Lemma~\ref{lem:basic}(b) gives $\prod_{l\ne k}f_l\in I$.
\end{proof}

The next result is purely algebraic; it holds in every commutative ring
and does not use $C(X)$.

\begin{proposition}
\label{prop:radical}
Let $I$ be a proper radical ideal of a commutative ring $R$ and
$n\ge1$. The following are equivalent.
\begin{enumerate}[label=\textup{(\alph*)}]
\item $I$ is $n$-absorbing.
\item $|\Min(I)|\le n$.
\item $I=P_1\cap\dots\cap P_m$ for some $m\le n$ and some prime ideals
$P_1,\dots,P_m$.
\end{enumerate}
If these hold and $I$ is written as in \textup{(c)} with
$P_1,\dots,P_m$ pairwise incomparable, then $\Min(I)=\{P_1,\dots,P_m\}$.
\end{proposition}

\begin{proof}
(c)$\Rightarrow$(a): Lemma~\ref{lem:elementary}(a), after repeating one
of the primes if $m<n$.

(b)$\Rightarrow$(c): $I$ is the intersection of the primes containing
it, and each of them contains a minimal one; hence
$I=\bigcap\Min(I)$, and $\Min(I)\ne\varnothing$.

(a)$\Rightarrow$(b): Suppose, to the contrary, that $P_1,\dots,P_{n+1}$
are distinct primes minimal over $I$. They are pairwise incomparable.
Let $S=R\setminus(P_1\cup\dots\cup P_{n+1})$, a multiplicative set
disjoint from $I$, and put $A=S^{-1}R$, $J=S^{-1}I$ and
$Q_i=S^{-1}P_i$. Then $J$ is a proper radical ideal of $A$, because
$\sqrt{S^{-1}I}=S^{-1}\sqrt I$.

\emph{Claim 1: $J$ is $n$-absorbing.} Let
$\alpha_i=a_i/s_i$ ($1\le i\le n+1$) satisfy $\prod\alpha_i\in J$. Then
$c\prod a_i\in I$ for some $c\in S$. Applying the $n$-absorbing property
of $I$ to the $n+1$ elements $ca_1,a_2,\dots,a_{n+1}$, we find $k$ with
$\prod_{l\ne k}a_l\in I$ (if $ca_1$ is omitted, $k=1$) or
$c\prod_{l\ne k}a_l\in I$ (if $a_k$ with $k\ge2$ is omitted). In both
cases $\prod_{l\ne k}\alpha_l\in S^{-1}I=J$, because $c$ becomes a unit.

\emph{Claim 2: $Q_1,\dots,Q_{n+1}$ are the only primes of $A$
containing $J$, they are pairwise incomparable, and
$J=Q_1\cap\dots\cap Q_{n+1}$.} A prime of $A$ containing $J$ has the
form $S^{-1}Q$ with $Q$ a prime of $R$, $I\subseteq Q$ and $Q\cap S=\varnothing$,
i.e.\ $Q\subseteq P_1\cup\dots\cup P_{n+1}$. By prime avoidance
$Q\subseteq P_i$ for some $i$, so $Q=P_i$ by minimality of $P_i$ over
$I$. Contraction gives back $P_i$ from $Q_i$, so the $Q_i$ are pairwise
incomparable. Finally $J$ is radical, hence the intersection of the
primes containing it.

Now, by prime avoidance and incomparability, choose
$b_i\in Q_i\setminus\bigcup_{j\ne i}Q_j$ for $i=1,\dots,n+1$. Then
$b_1\cdots b_{n+1}\in\bigcap_iQ_i=J$. By Claim~1 there is $k$ with
$\prod_{l\ne k}b_l\in J\subseteq Q_k$. As $Q_k$ is prime, some $b_l$ with
$l\ne k$ lies in $Q_k$, contradicting the choice of $b_l$.

Finally, if $P_1,\dots,P_m$ are pairwise incomparable primes with
$I=\bigcap P_i$, a prime $Q\supseteq I$ contains $\prod P_i$ and so
contains some $P_i$; hence the minimal primes over $I$ are exactly the
$P_i$.
\end{proof}

\begin{theorem}
\label{thm:main-equivalence}
Let $n\ge1$ and let $I$ be a proper $z$-ideal of $C(X)$. The following
are equivalent.
\begin{enumerate}[label=\textup{(\alph*)}]
\item $I$ is $n$-absorbing.
\item $\CZ[I]$ is an $n$-absorbing filter of $\CZ(X)$.
\item $|\Min(I)|\le n$.
\item $I=P_1\cap\dots\cap P_m$ with $m\le n$ and
$P_1,\dots,P_m$ pairwise incomparable prime ideals.
\end{enumerate}
Moreover, the set $\{P_1,\dots,P_m\}$ in \textup{(d)} is unique, equals
$\Min(I)$, and consists of $z$-ideals; and
$\CZ[I]=\CZ[P_1]\cap\dots\cap\CZ[P_m]$ with each $\CZ[P_i]$ a prime filter.
\end{theorem}

\begin{proof}
Since $I$ is a $z$-ideal it is radical (Lemma~\ref{lem:basic}(a)).
(a)$\Leftrightarrow$(b) is Proposition~\ref{prop:zero-set}; (a),
(c) and (d) are equivalent by Proposition~\ref{prop:radical}, noting that
in (c)$\Rightarrow$(d) one may discard primes that contain another one
of the family without changing the intersection.
The primes in $\Min(I)$ are $z$-ideals by Lemma~\ref{lem:min-z}. If
$I=\bigcap_{i\le m}P_i$ with the $P_i$ $z$-ideals and $A=\Z(f)\in\bigcap_i\CZ[P_i]$,
then $f\in P_i$ for all $i$ by Lemma~\ref{lem:basic}(b), so $f\in I$.
Hence $\bigcap\CZ[P_i]\subseteq\CZ[I]$; the other inclusion is clear,
and primality of $\CZ[P_i]$ is Lemma~\ref{lem:basic}(c).
\end{proof}

\begin{remark}
\label{rem:what-is-needed}
Statements (a), (c), (d) are equivalent for every radical ideal of every
commutative ring (Proposition~\ref{prop:radical}). The $z$-ideal
hypothesis is used for the zero-set formulation (b) and for the fact that
the primes $P_i$ are $z$-ideals. For $n=2$, the implication
(a)$\Rightarrow$(c) is already known in the stronger form that every
$n$-absorbing ideal of an arbitrary commutative ring has at most $n$
minimal primes~\cite[Theorem~2.5]{AndersonBadawi}. For $n=2$, see also
Badawi~\cite[Theorems~2.3 and~2.4]{Badawi}. The localization proof above
provides a self-contained alternative in the radical case.
\end{remark}

\section{Common zero sets}

\begin{lemma}
\label{lem:prime-K}
Let $P$ be a prime ideal of $C(X)$. Then $|\K(P)|\le1$. For $x\in X$ we
have $x\in\K(P)\iff P\subseteq M_x$, so $P\subseteq M_x$ iff
$\K(P)=\{x\}$. If $X$ is compact, then $\K(P)$ is a singleton.
\end{lemma}

\begin{proof}
Suppose $x\ne y$ both lie in $\K(P)$. Choose disjoint open sets
$U\ni x$, $V\ni y$ and $u,v\in C(X)$ with $u(x)=1$, $u=0$ off $U$,
$v(y)=1$, $v=0$ off $V$. Then $uv=0\in P$, so $u\in P$ or $v\in P$,
contradicting $x,y\in\K(P)$. The equivalence $x\in\K(P)\iff P\subseteq M_x$
is the definition. The last claim follows from the fact that proper
ideals of $C(X)$ are fixed when $X$ is compact.
\end{proof}

\begin{theorem}
\label{thm:K-bound}
Let $X$ be a Tychonoff space and $I$ an $n$-absorbing ideal of $C(X)$.
Then $|\K(I)|\le n$. The bound is sharp whenever $X$ has at least $n$ points: for distinct
$x_1,\dots,x_n\in X$, the ideal $M_{x_1}\cap\dots\cap M_{x_n}$ is an
$n$-absorbing $z$-ideal with common zero set $\{x_1,\dots,x_n\}$.
\end{theorem}

\begin{proof}
Suppose $x_1,\dots,x_{n+1}\in\K(I)$ are distinct. Choose pairwise
disjoint open sets $U_i\ni x_i$ and $u_i\in C(X)$ with $u_i(x_i)=1$ and
$u_i=0$ off $U_i$; thus at most one $u_i$ is nonzero at any point.
Put $f_j=\sum_{i\ne j}u_i$ for $j=1,\dots,n+1$.
At a point $p$, if all $u_i(p)=0$ then all $f_j(p)=0$, and if
$u_k(p)\ne0$ then $f_k(p)=0$. Hence $f_1f_2\cdots f_{n+1}=0\in I$.
Since $I$ is $n$-absorbing, $\prod_{j\ne m}f_j\in I$ for some $m$.
But for $j\neq m$ we have $f_j(x_m)=u_m(x_m)=1$, so this product equals
$1$ at $x_m$, whereas every element of $I$ vanishes at $x_m\in\K(I)$.
This contradiction proves $|\K(I)|\le n$.
Sharpness follows from Lemma~\ref{lem:elementary}(a), because the $M_{x_i}$
are prime $z$-ideals.
\end{proof}

\begin{proposition}
\label{prop:K-intersection}
For ideals $I_1,\dots,I_m$ of $C(X)$ we have $\K(I_1\cap\dots\cap I_m)=\K(I_1)\cup\dots\cup\K(I_m)$.
For every ideal $I$ we have $\K(\sqrt I)=\K(I)$.
\end{proposition}

\begin{proof}
The inclusion $\supseteq$ is clear. If $t\notin\bigcup_i\K(I_i)$, choose
$f_i\in I_i$ with $f_i(t)\ne0$; then $f_1\cdots f_m\in\bigcap_iI_i$ does not
vanish at $t$, so $t\notin\K(\bigcap_iI_i)$. For the second claim, $I\subseteq\sqrt I$ gives
$\K(\sqrt I)\subseteq\K(I)$; if $f\in\sqrt I$ then $f^k\in I$ vanishes
on $\K(I)$, hence so does $f$.
\end{proof}

\begin{proposition}
\label{prop:K-vs-Min}
Let $I$ be a proper $n$-absorbing $z$-ideal of $C(X)$ and write
$I=P_1\cap\dots\cap P_r$ with $\Min(I)=\{P_1,\dots,P_r\}$ and the $P_i$
pairwise incomparable prime $z$-ideals ($r\le n$). Then
\[
|\K(I)|=\bigl|\{\K(P_1),\dots,\K(P_r)\}\setminus\{\varnothing\}\bigr|,
\]
i.e.\ $|\K(I)|$ equals the number of distinct nonempty sets $\K(P_i)$.
\end{proposition}

\begin{proof}
By Proposition~\ref{prop:K-intersection},
$\K(I)=\bigcup_i\K(P_i)$. Each $\K(P_i)$ is either empty or a singleton
(Lemma~\ref{lem:prime-K}). Two singletons $\{x\}$ and $\{y\}$ coincide
iff $x=y$. Hence $|\K(I)|$ equals the number of distinct nonempty sets
$\K(P_i)$.
\end{proof}

\begin{remark}
For $n=2$, Theorem~\ref{thm:K-bound} also follows from Badawi's
description of $\sqrt I$ together with Lemma~\ref{lem:prime-K} and
Proposition~\ref{prop:K-intersection}. The direct proof above has the
advantage of being valid for every $n$ and of not using minimal primes.
\end{remark}

\section{On \texorpdfstring{$n$}{n}-absorbing ideals that are not \texorpdfstring{$z$}{z}-ideals}
\label{sec:nonz}

The preceding sections studied $n$-absorbing \emph{$z$}-ideals.
We now turn to arbitrary (not necessarily $z$-) $n$-absorbing ideals of
$C(X)$. We determine the number of minimal primes and characterize the
spaces in which every $n$-absorbing ideal is already a $z$-ideal. We also
prove that $z$-closure preserves $n$-absorption by passing to the
prime $z$-closures of the finitely many minimal primes.
A cautionary example is included to indicate the failure of naive
constructions; the remaining structural questions are collected as
open problems in Section~\ref{sec:open}.

\subsection{Minimal primes and the radical}

\begin{theorem}[Anderson--Badawi]
\label{thm:radical-nabs}
Let $n\ge1$ and let $I$ be a proper $n$-absorbing ideal of an arbitrary
commutative ring $R$. Then $|\operatorname{Min}(I)|\le n$. In particular
this holds for every $n$-absorbing ideal of $C(X)$, and $\sqrt{I}$ is the
intersection of at most $n$ prime ideals.
\end{theorem}

\begin{proof}
We first show that the radical $\sqrt I$ is $n$-absorbing. Let
$a_1\cdots a_{n+1}\in\sqrt I$. Then
$(a_1\cdots a_{n+1})^k=a_1^k\cdots a_{n+1}^k\in I$ for some $k\ge1$.
Since $I$ is $n$-absorbing, the product of some $n$ of the elements
$a_i^k$ lies in $I$; say $\prod_{l\ne m}a_l^k\in I$. Then
$\prod_{l\ne m}a_l\in\sqrt I$, so $\sqrt I$ is $n$-absorbing.

The ideal $\sqrt I$ is radical, satisfies $\Min(\sqrt I)=\Min(I)$, and
equals the intersection of its minimal primes. Applying
Proposition~\ref{prop:radical} to the radical $n$-absorbing ideal
$\sqrt I$ yields $|\Min(I)|=|\Min(\sqrt I)|\le n$, and
$\sqrt I=\bigcap_{P\in\Min(I)}P$.
\end{proof}

\subsection{The \texorpdfstring{$z$}{z}-closure of an \texorpdfstring{$n$}{n}-absorbing ideal}

For a prime ideal $P$ of $C(X)$, its $z$-closure
\[
P_z=\{f\in C(X):\text{ there is }g\in P
                  \text{ with }\Z(g)\subseteq\Z(f)\}
\]
is again a prime ideal. We include the argument because it is the
key point in proving that $z$-closure preserves $n$-absorption.

\begin{lemma}
\label{lem:prime-z-closure}
If $P$ is a prime ideal of $C(X)$, then $P_z$ is a proper prime
$z$-ideal.
\end{lemma}

\begin{proof}
Put
\[
\CZ[P]=\{\Z(f):f\in P\}.
\]
We first verify that $\CZ[P]$ is a proper filter of zero sets.
It contains $X=\Z(0)$, whereas $\varnothing\notin\CZ[P]$, since
a function with empty zero set is a unit and a proper ideal contains
no unit. If $\Z(f),\Z(g)\in\CZ[P]$, then
\[
\Z(f)\cap\Z(g)=\Z(f^2+g^2)\in\CZ[P],
\]
because $f^2+g^2\in P$. Moreover, if
$\Z(f)\in\CZ[P]$ and $\Z(f)\subseteq\Z(u)$, then
$fu\in P$ and
\[
\Z(fu)=\Z(f)\cup\Z(u)=\Z(u);
\]
hence $\Z(u)\in\CZ[P]$. Thus $\CZ[P]$ is upward closed.

We next prove that this filter is prime. Suppose
$\Z(f)\cup\Z(g)\in\CZ[P]$, and choose $h\in P$ such that
\[
\Z(h)=\Z(f)\cup\Z(g).
\]
Define
\[
a=(|f|-|g|)^+,\qquad b=(|g|-|f|)^+.
\]
Then $ab=0\in P$, so $a\in P$ or $b\in P$. Suppose that $a\in P$.
The function $h^2+a^2$ belongs to $P$, and
\[
\Z(h^2+a^2)=\Z(h)\cap\Z(a)\subseteq\Z(f).
\]
Indeed, if $t\in\Z(h)\cap\Z(a)$, then either $f(t)=0$, or
$g(t)=0$; in the latter case $a(t)=0$ gives
$|f(t)|\le |g(t)|=0$. By the upward-closure property just proved,
$\Z(f)\in\CZ[P]$. The case $b\in P$ is symmetric and gives
$\Z(g)\in\CZ[P]$. Thus $\CZ[P]$ is a prime filter.

By definition,
\[
P_z=\{u\in C(X):\text{ there is }v\in P
                    \text{ with }\Z(v)\subseteq\Z(u)\}.
\]
This is a $z$-ideal, and it is proper: if $1\in P_z$, some
$v\in P$ would satisfy $\Z(v)=\varnothing$, making $v$ a unit.
Furthermore,
\[
\CZ[P_z]=\CZ[P].
\]
Indeed, $P\subseteq P_z$ gives one inclusion, while the other follows
from the upward-closure property of $\CZ[P]$. Consequently,
$\CZ[P_z]$ is a prime filter. Lemma~\ref{lem:basic}(c), applied to
the $z$-ideal $P_z$, now shows that $P_z$ is prime.
\end{proof}

\begin{theorem}
\label{thm:z-closure-nabs}
Let $I$ be a proper $n$-absorbing ideal of $C(X)$, and write
\[
\Min(I)=\{P_1,\ldots,P_m\},\qquad m\le n.
\]
Then
\[
I_z=(P_1)_z\cap\cdots\cap(P_m)_z.
\]
Consequently, $I_z$ is an $n$-absorbing $z$-ideal and
\[
\K(I_z)=\K(I).
\]
\end{theorem}

\begin{proof}
Since $I$ is $n$-absorbing, its radical is $n$-absorbing and
\[
\sqrt I=P_1\cap\cdots\cap P_m
\]
by Theorem~\ref{thm:radical-nabs}. Each $(P_i)_z$ is a proper prime
$z$-ideal by Lemma~\ref{lem:prime-z-closure}.

If $f\in I_z$, choose $g\in I$ with $\Z(g)\subseteq\Z(f)$.
Since $I\subseteq P_i$, we have $g\in P_i$, and hence
$f\in(P_i)_z$ for every $i$. This proves
\[
I_z\subseteq\bigcap_{i=1}^m(P_i)_z.
\]

Conversely, suppose $f\in\bigcap_{i=1}^m(P_i)_z$. For each $i$,
choose $g_i\in P_i$ such that $\Z(g_i)\subseteq\Z(f)$, and put
$h_i=fg_i$. Then $h_i\in P_i$ and
\[
\Z(h_i)=\Z(f)\cup\Z(g_i)=\Z(f).
\]
Let $h=h_1\cdots h_m$. For each $i$, the factor $h_i$ belongs to
$P_i$, so $h\in\bigcap_iP_i=\sqrt I$. Thus $h^k\in I$ for some
positive integer $k$. Moreover,
\[
\Z(h^k)=\Z(h)=\bigcup_{i=1}^m\Z(h_i)=\Z(f),
\]
and therefore $f\in I_z$. This proves the reverse inclusion.

The displayed representation expresses $I_z$ as an intersection
of $m\le n$ prime ideals. Hence $I_z$ is $n$-absorbing; it is a
$z$-ideal because each $(P_i)_z$ is a $z$-ideal.

Finally, $I\subseteq I_z$ implies
$\K(I_z)\subseteq\K(I)$. If $f\in I_z$, choose $g\in I$ with
$\Z(g)\subseteq\Z(f)$. Then
\[
\K(I)\subseteq\Z(g)\subseteq\Z(f).
\]
Taking the intersection over all $f\in I_z$ gives
$\K(I)\subseteq\K(I_z)$, and equality follows.
\end{proof}

\begin{remark}
A direct attempt to prove that $I_z$ is $n$-absorbing by replacing
each factor from $I_z$ with a witness from $I$ does not suffice:
the resulting product of witnesses need not be related to the
hypothesis that the original product belongs to $I_z$. The proof
above instead uses the minimal primes of $I$ and the fact that the
$z$-closure of a prime ideal of $C(X)$ is prime.
\end{remark}

\subsection{When every \texorpdfstring{$n$}{n}-absorbing ideal is a \texorpdfstring{$z$}{z}-ideal}

\begin{corollary}[Classical $P$-space characterization]
\label{thm:every-nabs-is-z}
For a Tychonoff space $X$ and any fixed $n\ge1$ the following are
equivalent:
\begin{enumerate}[label=\textup{(\roman*)}]
\item $X$ is a $P$-space;
\item every ideal of $C(X)$ is a $z$-ideal;
\item every prime ideal of $C(X)$ is a $z$-ideal;
\item every $n$-absorbing ideal of $C(X)$ is a $z$-ideal.
\end{enumerate}
\end{corollary}

\begin{proof}
The equivalence of (i), (ii) and (iii) is classical; see
Gillman and Henriksen~\cite{GillmanHenriksen},
Gillman and Jerison~\cite[4J]{GillmanJerison}, and
Mason~\cite[Theorem~1.2]{Mason}. The implication (ii)$\Rightarrow$(iv) is
immediate since (iv) is a special case of (ii). Finally,
(iv)$\Rightarrow$(iii) because every prime ideal is $n$-absorbing.
\end{proof}

\subsection{A cautionary example}

A natural attempt to produce a non-$z$ $n$-absorbing ideal is to
intersect a prime ideal with a principal ideal. The following example
shows that such an ideal need not be $n$-absorbing at all, so this
construction cannot be used naively.

\begin{example}
\label{ex:principal-not-absorbing}
Let $X=[0,1]$ and let $I=(x)\subseteq C([0,1])$ be the principal ideal
generated by the identity function. Then $I$ is not $n$-absorbing for
any $n\ge1$.

For $n\ge2$: put $a_1=\dots=a_{n+1}=x^{1/(n+1)}$, the continuous real
$(n+1)$-st root of $x\ge0$. Then $a_1\cdots a_{n+1}=x\in I$, but the
$n$-fold product is $x^{n/(n+1)}\notin I$, because
$x^{n/(n+1)}/x=x^{-1/(n+1)}$ is unbounded at $0$ and hence not
continuous. For $n=1$: $I$ is not prime, since $x^{1/2}\cdot x^{1/2}=x\in I$
while $x^{1/2}\notin I$.
\end{example}

In particular, taking $P=M_0$ and $f=x$ in the example above gives
$P\cap(f)=(x)$, which is not $n$-absorbing. This shows that statements
such as ``$P\cap(f)$ is $n$-absorbing for every prime $P$ and every $f$''
are false in general: the intersection must be analyzed with respect to
the zero-set conditions satisfied by $f$ relative to $P$, which we leave
as an open problem in Section~\ref{sec:open}.

\subsection{Principal \texorpdfstring{$n$}{n}-absorbing ideals}

The preceding example shows that principal ideals need not be
$n$-absorbing. The next proposition gives an exact criterion for when
a principal ideal has this property.

\begin{proposition}
\label{prop:principal}
Let $f\in C(X)$ be a non-unit and let $n\ge1$. The following are
equivalent:
\begin{enumerate}[label=\textup{(\roman*)}]
\item $(f)$ is $n$-absorbing;
\item $\Z(f)$ is clopen and $|\Min((f))|\le n$.
\end{enumerate}
In particular, every principal $n$-absorbing ideal of $C(X)$ is a
$z$-ideal.
\end{proposition}

\begin{proof}
We first show that if $\Z(f)$ is not open, then $(f)$ is not
$n$-absorbing. Define continuous functions on $\mathbb R$ by
\[
r(t)=\operatorname{sgn}(t)|t|^{1/(n+1)}
\quad\text{and}\quad
s(t)=|t|^{1/(n+1)},
\]
and put $a_1=r\circ f$ and $a_2=\cdots=a_{n+1}=s\circ f$.
Then
\[
a_1a_2\cdots a_{n+1}=f\in(f).
\]
Every product $p$ of $n$ of these factors satisfies
\[
|p|=|f|^{n/(n+1)}.
\]
Suppose that such a product belongs to $(f)$, say $p=fg$ for some
$g\in C(X)$. On $\coz(f)$ we then have
\[
|g|=|f|^{-1/(n+1)}.
\]
Choose $x\in\Z(f)\setminus\operatorname{int}\Z(f)$. For every
$M>0$, continuity of $f$ gives a neighbourhood $U$ of $x$ on which
$|f|<M^{-(n+1)}$. Since $x\notin\operatorname{int}\Z(f)$,
the set $U\cap\coz(f)$ is nonempty. At every point of this set,
$|g|>M$. Thus $g$ is unbounded in every neighbourhood of $x$,
contrary to continuity at $x$. Hence no product of $n$ of the
$a_i$ belongs to $(f)$, and $(f)$ is not $n$-absorbing.

Now suppose that $(f)$ is $n$-absorbing. The first part shows that
$\Z(f)$ is open; since it is a zero set, it is also closed. Let
$e=\chi_{\coz(f)}$, which is continuous because $\coz(f)$ is clopen.
The function $u=f+1-e$ has no zeros, so it is a unit, and $f=eu$.
Consequently
\[
(f)=(e)=\{g\in C(X):g=ge\}
     =\{g\in C(X):\Z(g)\supseteq\Z(f)\}.
\]
Thus $(f)$ is a $z$-ideal. Theorem~\ref{thm:main-equivalence}
now gives $|\Min((f))|\le n$.

Conversely, assume that $\Z(f)$ is clopen and
$|\Min((f))|\le n$. The same characteristic-function argument shows
that $(f)$ is a proper $z$-ideal. Theorem~\ref{thm:main-equivalence}
then implies that $(f)$ is $n$-absorbing.
\end{proof}

\section{Open problems}
\label{sec:open}

\begin{problem}
\label{prob:classify-gap}
Is every $n$-absorbing non-$z$ ideal of $C(X)$ of the form
\[
Q\cap(f_1)\cap\dots\cap(f_k),
\]
where $Q$ is an $n$-absorbing $z$-ideal and $f_1,\dots,f_k\in C(X)$?
If not, what is the simplest counter-example? (Example~\ref{ex:principal-not-absorbing}
shows that the answer is not trivially affirmative.)
\end{problem}

\begin{problem}
\label{prob:exact-generators}
Let $Q=P_1\cap\dots\cap P_m$ be an $n$-absorbing $z$-ideal. Give
necessary and sufficient zero-set conditions on elements
$f_1,\dots,f_k\in C(X)$ such that
$I=Q\cap(f_1)\cap\dots\cap(f_k)$ satisfies $I_z=Q$.
\end{problem}

\begin{problem}
\label{prob:exact-condition}
Find a purely topological (zero-set) condition on the family
$\{f_1,\dots,f_k\}$ that is necessary and sufficient for
\[
\bigl(Q\cap(f_1)\cap\dots\cap(f_k)\bigr)_z=Q
\]
to hold, without assuming first-countability.
\end{problem}

\end{document}